\documentclass[a4paper,11pt]{article}
\usepackage{charter}
\usepackage{latexsym}
\usepackage{amsmath}
\usepackage{amsfonts}
\usepackage{amssymb}
\usepackage{vmargin}
\usepackage{url}
\usepackage{undertilde}
\usepackage{color}
\usepackage{booktabs}
\usepackage[singlelinecheck=false]{caption}
\setmarginsrb{2cm}{2cm}{2cm}{2cm}{0mm}{0mm}{0mm}{10mm}

\begin{document}
\title{\vspace*{-0.7cm}\textbf{Some pros and cons of implementing parallel and block teachings for mathematics modules}}
\author{\normalsize N. Karjanto{\small $^{1}$} and S. T. Yong{\small $^{2}$}\\
\textsl{\small $^{1}$Department of Applied Mathematics \; $^{2}$Department of Foundation in Engineering}\\
\textsl{\small Faculty of Engineering, The University of Nottingham Malaysia Campus}\\
\textsl{\small Jalan Broga, Semenyih 43500, Selangor, Malaysia}}
\date{\footnotesize }
\maketitle

\begin{abstract}
\noindent
The Department of Applied Mathematics at the University of Nottingham Malaysia Campus has a responsibility for delivering mathematics modules for engineering students. Due to the significantly large number of students, two methods of teaching delivery--parallel teaching and block teaching--have been implemented. This article discusses some pros and cons between these two methods, particularly for the Foundation programme and the first year of the Undergraduate programme in Engineering. Whether parallel teaching or block teaching is implemented, feedback comments from the students indicate that some areas need to be paid attention to. \\

\setlength{\parindent}{0pt}
\end{abstract}

\section{Introduction}

This article reveals some pros and cons when implementing different methods of teaching delivery in several undergraduate and foundation mathematics modules at our university. Our campus houses four different faculties and our department belongs to the Faculty of Engineering. Although the Department of Applied Mathematics does not really offer any study programme in math, it plays an important role as a service department not only for the Undergraduate programme but also for the Foundation programme in Engineering. Currently, there are four engineering departments housed in the Faculty of Engineering: Electrical and Electronic Engineering, Chemical and Environmental Engineering, Mechanical, Materials and Manufacturing Engineering and Civil Engineering.

Since our university is a branch campus of the University of Nottingham system, we adopt the English higher education system in our curricula, where the undergraduate level is divided into three years~\cite{framework}, in contrast to the North American system where the students will normally spend four years in their undergraduate study. This means that the students who enter Year~1 have an entry requirement of A~(advanced)-level or a similar qualification. However, a significant number of our potential students posses only an O~(ordinary)-level qualification. In order to accommodate them, the faculty arranges a special programme known as the Foundation in Engineering. This programme lasts for three semesters which normally starts around April or July and ends in May of the following year. The purpose of this programme is to equip the students with merely an O-level qualification to enter the Undergraduate programme in Engineering. In this case, at the end of the third semester of their Foundation study, they would have a similar level to the students with an A-level qualification. After completing the Foundation programme and satisfying the progression rules, the students may proceed to choose one study programme offered by the four engineering departments.

Our department also services a number of higher level undergraduate and graduate mathematics modules, such as for Year 3 and Year 4 (MSc and MEng) students. These modules are usually offered as elective modules and generally, the number of students who take such optional modules is relatively small in comparison to the compulsory math modules. Keeping this in mind, the math modules discussed in this article will refer to the modules in both Foundation and Year 1 of the Undergraduate programmes.

This article is organized as follows. The following section will briefly mention the number of students we currently have and we expect to have in the future. It will also explain how these large numbers of students are organized into groups. Section~\ref{delivery} will discuss the difference between the two methods of teaching delivery, the parallel teaching and block teaching. In Section~\ref{procon}, we will reveal some pros and cons between parallel teaching and block teaching methods, both from lecturers' and students' perspectives. Finally, Section~\ref{conclude} draws a conclusion and provides remark to our discussion, with a balanced view regarding the two teaching methods.

\section{Class organization}

Due to a significantly large number of students and a limited number of academic staff in our department, methods of teaching delivery are of a great concern to both sides of the teaching activity, the students and the lecturers. Depending on the academic year, the number of students may vary from one cohort to the other. The number of students currently registered at the Foundation programme is around~250. However, we may anticipate this number to increase up to 300~students or more in the coming academic years. In Year~1 of the Undergraduate programme, the total number of students has grown to around~360, due to additional students registered straight away to Year~1 instead of passing through the Foundation stage. These newcomers normally have done their A-level qualification somewhere else. Similarly, we might also want to anticipate an increase in this number to around 400~students or more in the future. It should be taken into account that those who are taking mathematics modules in Year~1 may also include the students from another faculty, for instance from the Business School, whereby the students take these as their elective modules.

Regarding the class arrangement, the Year~1 students of the Undergraduate programme have been clustered into their own respective departments. For the students in the Foundation programme, however, there is not really a clear mechanism regarding the arrangement of the class. In this case, we adopt the grouping system that has been arranged by our colleagues at the School of English Language Education, Faculty of Arts and Education. Whatever grouping system we adopt, the same principle applies, i.e. the Foundation students are divided into three different groups with the same number of students. Thus, each group contains around 80~students. However, we are not able to maintain a similar proportion anymore for Year~1 students. This is caused by the fact that each department decides the number of students it is willing to admit and some trends that develop among the Foundation students regarding which department they would enrol. Thus, the four groups of Year~1 students could range from 50 to 140 students.

Until the recent academic year, it is compulsory for the students enrolled in the Foundation programme to take three mathematics modules. In Semester~1, the module is HG1BMT~Basic Mathematical Techniques (10 credit points). At this stage, the students are enrolled at two different intakes (April and July) and the number of students apparently is not so large. In Semester~2 and~3, the module is HG1FND~Foundation Mathematics (20 credit points), which is delivered as a year-long module, where each semester carries 10 credit points. In Semester~3, another math module is added, i.e. HG1M02~Applied Algebra (10 credit points). All Year~1 Engineering students must take HG1M11~Engineering Mathematics~1 and HG1M12~Engineering Mathematics~2 in the Autumn and Spring semesters, respectively. Both modules have the weight of 10 credit points each. It is essential to note that although the Year~1 Engineering students have chosen their own respective departments, the contents of the material in both modules mentioned above are identical.

In the following section, we discuss two methods of teaching delivery in connection to the mathematics modules mentioned above.

\section{Teaching delivery} \label{delivery}

According to~\cite{Friend93}, parallel and block teaching (also referred to as `station teaching') are two of the five types of co-teaching delivery for secondary school level. The other three are `lead and support', `alternative teaching' and `team teaching'.

In the context of this article, `parallel teaching' refers to the teaching delivery given by one lecturer handling one group of students only for the whole semester, or even the whole year. Although the lecturers jointly plan for the teaching material so that they will deliver the contents in a similar way, it is inevitable that some variations will occur from one group to the other since every lecturer has his/her own teaching style, method and approach.

On the other hand, the `block teaching' refers to the teaching delivery given by two or more different lecturers to each group of students. This means that within one semester, each lecturer will handle only certain topics that cover a couple weeks of material for the whole cohort of students. So in one semester the students will enjoy being taught by several lecturers. Some examples of block teaching plans are given in Appendix~A. In our example, three lecturers are involved in one particular module. Block teaching by two lecturers is common too, but we have never implemented four lecturers or more in the block teaching.

For each module that employs the block teaching plan, there is one lecturer assigned as a module convener. He/she will organize the other members of the group on who will teach which weeks. In our example, we run an 11-week semester. Since the final week is normally used for a test, then there will be only 10 weeks of teaching activities. In this case, one lecturer, normally the module convener, will handle four weeks of lectures, while the other two will handle three weeks of class. In addition, the module convener also distributes the workload to the other members of the group. For instance, it has to be determined who will contribute which problems for quizzes, tests or the final exam.

For every 10 credit points module, there are a two-hour lecturing session and a one-hour tutorial session. During the two-hour lecturing session, the lecturer covers the specified topics of material. Depending on the availability of the time, during this session, the lecturer might also employ students' participation, for example by encouraging them to solve some relevant problems. During the one-hour tutorial session, which is usually conducted on a different day than the lecture, the lecturer discusses the problem sheets given earlier to the students. During this session, the students have ample opportunity to practice their problem-solving skills. The lecturer gives a measure of guidance to the students when they attempt to solve more challenging problems and even might demonstrate a step by step explanation in solving a certain mathematical problem. It is good to note that the tutorial session is not necessarily conducted by the same lecturer as the one who conducts the lecturing session.

\section{Pros and cons of parallel and block teaching} \label{procon}

We have implemented both parallel and block teaching in delivering several mathematics modules during the past couple of years. In order to find out the general opinions, including some pros and cons between these two teaching methods, we have conducted a number of interview sessions with several of our lecturer colleagues and a number of students from the Foundation programme in Engineering.

After implementing parallel teaching delivery, we received a number of reports from the students regarding the different contents of material they received from one group to another. It is observed that the students tend to compare the teaching quality of their lecturers and the worked examples that they receive with their own friends from another class. When the students find out that they receive less worked examples or different ones, they might feel that they have been treated sufficiently unfairly that they have grounds for  complaint. With this in mind, in the current academic year, we implemented block teaching for all mathematics modules for both Foundation and Year~1 levels.

The following are some pros and cons regarding the two methods of teaching delivery. The parallel teaching delivery allows the lecturer to get to know the students better during the semester period. Although we always deal with a large number of students, knowing the ability of each student is particularly useful in smaller class size when guiding and helping them through their difficulties. From the lecturer's viewpoint, particularly the newly appointed staff who lack some teaching experience, parallel teaching is also useful as an opportunity to learn the module materials for the complete semesters, instead of merely a few topics only. Nevertheless, despite the efforts around jointly planning the content delivery in a similar way, some variations are simply inevitable, particularly when the lecturers have a wide range of teaching experience.

On the other hand, the block teaching allows identical deliveries to all  students regardless of the group in which they are sitting. By implementing block teaching, all groups of students will receive the same lecture notes, worked examples, example sheets and exercises and also will experience the same teaching style from a particular lecturer. The block teaching delivery exposes the students to many different lecturers in our department and allows them to taste the diversity of different teaching styles.

However, the block teaching delivery is not without shortcomings. We again receive some reports from the students mentioning that they find a difficulty in adjusting to different styles of teaching delivery from one lecturer to the another. In particular, the Foundation students suffer the most since the one-hour tutorial for discussing example sheets and exercises is usually being handled by a different lecturer.

It is common that different lecturers use different notation to express certain concepts or definitions. This is not surprising since when we check several mathematics textbooks, we will discover that different authors will use different mathematical symbols to express certain quantities. For example, many textbooks have used a standard notation for a vector using a boldface, i.e. $\mathbf{u}$. However, there are some variations to express a vector. For blackboard notation, it is acceptable to denote a vector using an arrow on top of the letter, i.e. $\vec{\textmd{\,u}}$. Some, however, use an underline or an undertilde notation to denote a vector, i.e. \underline{u} or $\utilde{\textmd{u}}$.

It has also been reported that some tutorial lecturers often have a completely different approach in solving certain problems compared to the main lecturer who is teaching them in the two-hour lecture session. As a consequence, many students get confused regarding which methods they should use in solving the problems in the exams. Some examples of these differences are given in Appendix~B. The difference in notations and techniques used across mathematics modules in the context of engineering lectures has also been observed in an Australian university~\cite{Varsavsky95}. The author gives an example that Cramer's rule appears to be widely used to solve systems of linear equation, but only Gaussian elimination and matrix inversion are taught in the mathematics subjects.

Another major concern when adopting the block teaching method is in the context of composing a timetable. One colleague who handles the timetable has said that block teaching poses more restrictions to certain hours. For instance, Foundation Mathematics module is scheduled at 9 am--11 am on Monday for Group~A, on Tuesday for Group~B and on Wednesday for Group~C. Block teaching will require two or three different time slots for different days since one lecturer may only handle one group at each time. On the contrary, parallel teaching only requires one teaching slot for the whole groups, let say on Monday 9 am--11 am, and the two or three lecturers simply go to their own respective groups. This timetabling problem can be solved by implementing a certain scheduling program. An example on building university timetabling using certain mathematical programming is given in~\cite{Gueret96}. For school timetable construction, the interested readers may consult~\cite{Willemen02}.

\section{Conclusion and Remark} \label{conclude}

In this article, we have discussed some pros and cons regarding two different methods of delivering mathematics modules to engineering students at the University of Nottingham Malaysia Campus. Implementing the parallel teaching method has given a good opportunity for both sides of lecturers and students in getting to know each other better so that the lecturers could monitor the progress many of their students, if not all, until the semester ends. On the other hand, parallel teaching gives room for more variation in teaching techniques and delivery and in turn, the students may feel that they are not treated sufficiently fairly since they do not receive the same materials as their counterparts in the other classes do.

With this in mind, the block teaching has been introduced in the past academic year. The number of complaints has been reduced after implementing this method. However, a number of shortcomings of block teaching has also been identified. Students find difficulties in adjusting the teaching style of one lecturer to the other within a short period of time, i.e. a couple of weeks period. Some students confuse with different notations being used within the same module and different techniques introduced by different lecturers to solve certain mathematical problems.

Whether we would like to implement parallel teaching or block teaching delivery in the future, it is important that all lecturers who are involved sit down and discuss the matter beforehand. Hopefully, this will prevent a significant variation in the teaching materials. Furthermore, in order to reduce students' confusion, special attention must be given to the notations and techniques used in the modules. Finally, by implementing a particular scheduling program, it will reduce the trouble that comes when arranging the timetable for block teaching delivery.

\section*{\normalsize Acknowledgement}
{\small We would like to thank Dr. Rohaizan Osman as the Director of Studies of Department of Applied Mathematics, Mr. Chiang Choon Lai as the timetable coordinator for Foundation Studies and a number of students at the Faculty of Engineering who are willing to express their opinions regarding the different methods of teaching delivery.}

\newpage
\subsection*{Appendix A \; Some examples of block teaching plan}
\begin{table}[htbp]
\flushleft
\caption{Block teaching plan for Foundation Mathematics (Autumn/Fall)}
\begin{tabular}{@{}cp{11cm}l@{}}
\toprule
Week & Topics & Lecturers \\ \hline 
  1  & Introduction to trigonometric function &  \\
  2  & Identities of trigonometric functions & Dr. X. Analyst \\
  3  & Solving trigonometric equations &  \\ \hline
  4  & Introduction to differentiation &  \\
  5  & Basic rules of differentiation & Dr. Y. Topologist \\
  6  & Product, quotient and chain rules &  \\ \hline
  7  & Derivatives of trigonometric function &  \\
  8  & Derivatives of exponential and logarithmic functions & Prof. Z. Mathematician\\
  9  & Derivatives of implicit functions & (Module convener) \\
 10  & Introduction to complex numbers  & \\ \hline
 11  & Test & \\ 
\bottomrule
\end{tabular}
\end{table}

\begin{table}[htbp]
\caption{Block teaching plan for Foundation Mathematics (Spring)}
\begin{tabular}{@{}cp{11cm}l@{}}
\toprule
Week & Topics & Lecturers \\ \hline
  1  & Applications of derivative 1 &  \\
  2  & Applications of derivative 2 & Dr. V. Algebraist \\
  3  & Curve sketching 1 &  \\ \hline
  4  & Curve sketching 2 &  \\
  5  & Introduction to integral & Prof. Z. Mathematician \\
  6  & Techniques of integration & (Module convener)  \\
  7  & Application of integral 1 &  \\ \hline
  8  & Application of integral 2 & \\
  9  & Numerical integration & Dr. W. Statistician \\
 10  & Binomial theorem for any rational index  & \\ \hline
 11  & Test & \\ 
\bottomrule
\end{tabular}
\end{table}

\begin{table}[h!]
\caption{Block teaching plan for Engineering Mathematics 1}
\begin{tabular}{@{}cp{11cm}l@{}}
\toprule
Week & Topics & Lecturers \\ \hline
  1  & Coordinates system &  \\
  2  & Algebra of complex numbers & Prof. G. Geometer  \\
  3  & Matrix algebra &  (Module convener) \\
  4  & System of equations and eigenvalue problems &  \\ \hline
  5  & Application of matrix algebra &  \\
  6  & Functions and their properties & Dr. N. Number-Theorist   \\
  7  & Advanced differential calculus &  \\ \hline
  8  & Taylor series & \\
  9  & Curve sketching & Dr. Q. Numerical-Analyst \\
 10  & Integral calculus of one variable  & \\ \hline
 11  & Test & \\ 
\bottomrule
\end{tabular}
\end{table}

\begin{table}[h!]
\begin{center}
\caption{Block teaching plan for Engineering Mathematics 2}
\begin{tabular}{@{}cp{10cm}l@{}}
\toprule
Week & Topics & Lecturers \\ \hline
   1 & Function of several variables &  \\
   2 & Partial derivatives and their interpretations & Dr. M. Mathematical-Physicist \\
   3 & Chain rule and total derivatives &   \\ \hline
   4 & First order ordinary differential equations &  \\
   5 & Exact differential equations &  Prof. G. Geometer \\
   6 & Boundary and initial conditions & (Module convener)   \\
   7 & Second-order linear constant coefficient & \\
     & ordinary differential equations & \\ \hline
   8 & Vector spaces & \\
   9 & Scalar and vector products & Dr. T. Graph-Theorist \\
  10 & Applications of vector methods in engineering & \\ \hline
  11 & Test & \\ 
\bottomrule  
\end{tabular}
\end{center}
\end{table}

\subsection*{Appendix B \; Some examples of problems solved by different methods}

\subsubsection*{Example 1 \; Vector (cross) product}

To evaluate the vector product (or vector product) between two vectors, we need to calculate the determinant of a 3 $\times$ 3 matrix. See~\cite{Anton05} or any other books on linear algebra for more complete descriptions on matrix determinant. There are at least two ways to calculate this determinant, i.e. using a combination of three 2 $\times$ 2 determinants and the other one is using Sarrus' rule or Sarrus' scheme. To see the difference between these methods, let us see an example. Let $\mathbf{a} = a_{1} \mathbf{i} + a_{2} \mathbf{j} + a_{3} \mathbf{k}$ and $\mathbf{b} = b_{1} \mathbf{i} + b_{2} \mathbf{j} + b_{3} \mathbf{k}$, then the cross product of $\mathbf{a}$ and $\mathbf{b}$ reads
\begin{equation*}
  \mathbf{a} \times \mathbf{b} = \left|%
\begin{array}{ccc}
  \mathbf{i} & \mathbf{j} & \mathbf{k} \\
  a_{1}   & a_{2}   & a_{3} \\
  b_{1}   & b_{2}   & b_{3} \\
\end{array}%
\right|.
\end{equation*}
Using the first method, the cross product is given by
\begin{align*}
\mathbf{a} \times \mathbf{b} &= \left|%
\begin{array}{cc}
  a_{2} & a_{3} \\
  b_{2} & b_{3} \\
\end{array}%
\right| \mathbf{i} - \left|%
\begin{array}{cc}
  a_{1} & a_{3} \\
  b_{1} & b_{3} \\
\end{array}%
\right| \mathbf{j} + \left|%
\begin{array}{cc}
  a_{1} & a_{2} \\
  b_{1} & b_{2} \\
\end{array}%
\right| \mathbf{k} \\
&= (a_{2}b_{3} - a_{3}b_{2}) \mathbf{i} - (a_{1}b_{3} - a_{3}b_{1})
\mathbf{j} + (a_{1}b_{2} - a_{2}b_{1}) \mathbf{k}.
\end{align*}

Sarrus' rule is a method of evaluating the determinant of a 3 $\times$ 3 matrix. This method is named after the 19$^{\textmd{th}}$-century French mathematician Pierre
Fr\'{e}d\'{e}ric Sarrus (1798-1861). The scheme is given as follows. Repeat the first two columns of the matrix behind the third column, so that we have five columns in a row. Then add the products of the diagonal components from top to bottom and subtract the products of the diagonal components going from bottom to top~\cite{wikiSarrus}. Thus, using Sarrus' rule, the cross product is given by
\begin{align*}
\mathbf{a} \times \mathbf{b} &= \left|
\begin{array}{cccccc}
  \mathbf{i} & \mathbf{j} & \mathbf{k} & \mid & \mathbf{i} & \mathbf{j} \\
  a_{1} & a_{2} & a_{3} & \mid & a_{1} & a_{2} \\
  b_{1} & b_{2} & b_{3} & \mid & b_{1} & b_{2} \\
\end{array}
\right| 
 = a_{2} b_{3} \mathbf{i} + a_{3} b_{1} \mathbf{j} + a_{1} b_{2} \mathbf{k} - b_{1} a_{2} \mathbf{k} - b_{2} a_{3} \mathbf{i} - b_{3} a_{1} \mathbf{j}\\
&= (a_{2}b_{3} - a_{3}b_{2}) \mathbf{i} + (a_{3}b_{1} - a_{1}b_{3})
\mathbf{j} + (a_{1}b_{2} - a_{2}b_{1}) \mathbf{k}.
\end{align*}
We can observe that in fact, both rules produce the same result. A confusion might occur at the final expressions. Using Sarrus' rule, all terms have positive signs while the earlier method places a negative sign in front of the second term.

\subsubsection*{Example 2 \; Definite integral by substitution}

Another example of implementing different methods is to evaluate a definite integral using substitution method. See~\cite{Smith08} or any other Calculus books for more detailed explanations and examples. Suppose we want to evaluate
\begin{equation*}
I = \int_{-\sqrt{\pi/6}}^{\sqrt{2\pi/3}} \; \frac{7}{5}x \cos \left(x^{2} + \frac{\pi}{6} \right) \, {\rm d}x.
\end{equation*}
Using the substitution method, we let $u = x^{2} + \pi/6$, so ${\rm d}u = 2x \, {\rm d}x$ or $x \, {\rm d}x = \frac{1}{2}{\rm d}u$. But we have also to take care the boundary of integration, i.e. for $x = -\sqrt{\pi/6}$, $u = \pi/3$ and for $x = 2\pi/3$, $u = 5\pi/6$. Thus the integral becomes
\begin{equation*}
  I = \int_{\pi/3}^{5\pi/6} \frac{7}{10} \cos u \, {\rm d}u = \frac{7}{10} \sin u \bigg|_{\pi/3}^{5\pi/6} = \frac{7}{20}(1 - \sqrt{3}).
\end{equation*}
A similar method is without introducing a new variable $u$ and therefore there is no need to change the boundaries of integration. We can write the definite integral above as follows:
\begin{equation*}
I = \int_{-\sqrt{\pi/6}}^{\sqrt{2\pi/3}} \; \frac{7}{10} \cos \left(x^{2} + \frac{\pi}{6} \right) \; {\rm d}\left(x^{2} + \frac{\pi}{6} \right) 
  = \left. \frac{7}{10} \sin \left(x^{2} + \frac{\pi}{6} \right) \right|_{-\sqrt{\pi/6}}^{\sqrt{2\pi/3}} 
  = \frac{7}{20}(1 - \sqrt{3}).
\end{equation*}
It is observed that the students have a better grasp on the method of substitution by introducing a new variable $u$. However, many of them often forget to replace the boundaries of integration in term of the new variable. On the other hand, the latter method is useful since the students do not have to worry about the boundaries of integration. However, rewriting the integral in the desired form as the above can be quite challenging to many students in general.

{\small
\begin{thebibliography}{99}
\bibitem{framework} The framework for higher education qualifications in England, Wales and Northern Ireland. Available online at \url{https://www.ncl.ac.uk/ltds/assets/documents/qsh-progapp-fheqsummary.pdf}. Accessed \today.

\bibitem{Friend93} Friend, M., Reising, M., \& Cook, L. (1993). Co-teaching: An overview of the past, a glimpse at the present, and considerations for the future. {\it Preventing School Failure}, {\bf 37}(4): 6--10.

\bibitem{Varsavsky95} Varsavsky, C. (1995). The design of the mathematics curriculum for engineers: A joint venture of the mathematics department and the engineering faculty. {\it European Journal of Engineering Education}, {\bf 20}(3): 341--345.

\bibitem{Gueret96} Gu\'{e}ret, C., Jussien, N., Boizumault, P. and Prins, C. (1996). Building university timetables using constraint logic programming. {\it Practice and Theory of Automated Timetabling}, Springer Berlin/Heidelberg.

\bibitem{Willemen02} Willemen, R. J. (2002). {\it School timetable construction--algorithms and complexity}, PhD Thesis, Technical University Eindhoven, The Netherlands.

\bibitem{Anton05} Anton, H. (2005). {\it Elementary Linear Algebra}, 9$^{\textmd{th}}$ edition, John Wiley \& Sons, Inc. 

\bibitem{wikiSarrus} \url{http://en.wikipedia.org/wiki/Rule_of_Sarrus}. Accessed \today.

\bibitem{Smith08} Smith, R. T. and Minton, R. B. (2008). {\it Calculus}, 3$^{\textmd{rd}}$ edition, McGraw-Hill.
\end{thebibliography}
}
\end{document}